\documentclass[a4paper]{amsart}
\def\myfbox#1{#1}

\begin{document}
\title[Prime knots up to arc index 11]%
      {A tabulation of prime knots up to arc index 11}
\author[Gyo Taek Jin and Wang Keun Park]{Gyo Taek Jin and Wang Keun Park}
\thanks{}
\address{Department of Mathematical Sciences, KAIST, Daejeon 305-701 Korea}
\email{trefoil@kaist.ac.kr, lg9004@hotmail.com}
\begin{abstract}
As a supplement to the authors' article~\cite{Jin2007}, we present minimal arc presentations of the prime knots up to arc index 11.
\end{abstract}
\maketitle

\section{Introduction}
An \emph{arc presentation\/} of a knot is an embedding of a knot into the union of finitely many vertical half planes whose common boundary is the $z$-axis so that each of the half planes intersect the knot in a single properly embedded curve. These curves are called the \emph{arcs\/} of the arc presentation. The minimal number of arcs among all arc presentations of a given knot is called the \emph{arc index.}

\begin{figure}[h]
\centering
\setlength{\unitlength}{0.6mm}
\begin{picture}(40,40)
\thicklines \put(0,0){\line(1,0){30}} \put(30,0){\line(0,1){30}}
\put(30,30){\line(-1,0){7}} \put(17,30){\line(-1,0){7}}
\put(10,30){\line(0,-1){20}} \put(10,10){\line(1,0){17}}
\put(33,10){\line(1,0){7}} \put(40,10){\line(0,1){30}}
\put(40,40){\line(-1,0){20}} \put(20,40){\line(0,-1){20}}
\put(20,20){\line(-1,0){7}} \put(7,20){\line(-1,0){7}}
\put(0,20){\line(0,-1){20}}
\end{picture}
\caption{A grid diagram of a trefoil knot}
\end{figure}

\begin{figure}[h]
\centering
\setlength{\unitlength}{0.7mm}
\begin{picture}(50,55)(-10,0)
{
\put(-10,10){\line(1,-1){10}}
\put(-10,10){\line(4,-1){7.5}} \put(2.5,6.875){\line(4,-1){27.5}}
\put(-10,20){\line(2,-1){7.5}} \put(2.5,13.75){\line(2,-1){7.5}}
\put(-10,20){\line(5,-1){7.5}} \qbezier(2.5,17.6)(5,17.2)(7.5,16.6)
                               \put(12.5,15.4){\line(5,-1){15}} \put(32.5,11.6){\line(5,-1){7.5}}
\put(-10,30){\line(1,-1){10}}
\put(-10,30){\line(3,-1){17.5}} \put(12.5,22.5){\line(3,-1){7.5}}
\put(-10,40){\line(2,-1){20}}
\put(-10,40){\line(4,-1){27.5}} \put(22.5,31.875){\line(4,-1){7.5}}
\put(-10,50){\line(5,-1){50}}
\put(-10,50){\line(3,-1){30}}
\put(30,0){\line(0,1){30}}
\put(10,30){\line(0,-1){20}}
\put(40,10){\line(0,1){30}}
\put(20,40){\line(0,-1){20}}
\put(0,20){\line(0,-1){20}}
}
\qbezier[50](-10,5)(-10,30)(-10,55)
\linethickness{0.3mm}
\put(0.7,0){\line(1,0){28.6}}
\put(29.3,30){\line(-1,0){6.3}}
\put(17,30){\line(-1,0){6.3}}
\put(10.7,10){\line(1,0){16.3}}
\put(33,10){\line(1,0){6.3}}
\put(39.3,40){\line(-1,0){18.6}}
\put(20.7,20){\line(-1,0){7.7}}
\put(7,20){\line(-1,0){6.3}}
\put(29.3,0){\line(0,1){30}}
\put(10.7,30){\line(0,-1){20}}
\put(39.3,10){\line(0,1){30}}
\put(20.7,40){\line(0,-1){20}}
\put(0.7,20){\line(0,-1){20}}
\end{picture}
\caption{Construction of an arc presentation from a grid diagram}\label{fig:grid2arc}
\end{figure}

A \emph{grid diagram\/} is a knot diagram whose projection is a closed curve which is composed of finitely many horizontal line segments and the same number of vertical line segments such that vertical line segments always cross over horizontal line segments.
A grid diagram can be converted easily to an arc presentation with the number of arcs equal to the number of vertical line segments as indicated in Figure~\ref{fig:grid2arc}.
Conversely, an arc presentation can be easily converted to a grid diagram with the number of vertical line segments equal to the number of arcs.

Table~1 shows the number of prime knots of given arc index and minimal crossing number. Though we tried our best to make the table as accurate as possible, some of the knots counted by the four italicized numbers, 75, 12, 3, 3, may have smaller crossing numbers but not smaller than 17.

{\small
\begin{table}[h]
\newcommand{\vsp}{\vrule width0pt height 10pt depth.5pt}
\begin{tabular}{|c|c|c|c|c|c|c|c|}
\hline
\vrule width0pt height 12pt depth22pt
\lower15pt\hbox{Crossings} \kern5pt Arc index &
        \hbox to 20pt{\hfil 5\hfil}&
        \hbox to 20pt{\hfil 6\hfil}&
        \hbox to 20pt{\hfil 7\hfil}&
        \hbox to 20pt{\hfil 8\hfil}&
        \hbox to 20pt{\hfil 9\hfil}&
        \hbox to 20pt{\hfil {10}\hfil}&
        \hbox to 20pt{\hfil {11}\hfil}\\
\hline
\vsp  3&1& & & & &  &  \\
\hline
\vsp  4& &1& & & &  &  \\
\hline
\vsp  5& & &2& & &  &  \\
\hline
\vsp  6& & & &3& &  &  \\
\hline
\vsp  7& & & & &7&  &  \\
\hline
\vsp  8& & &1&2& &{18}&  \\
\hline
\vsp  9& & & &2&6&  &{41}\\
\hline
\vsp 10& & & &1&9&32&  \\
\hline
\vsp 11& & & & &4&46&135 \\
\hline
\vsp 12& & & & &2&48& 211\\
\hline
\vsp 13& & & & & &49& 399  \\
\cline{1-8}
\vsp 14& & & & & &17& 477  \\
\cline{1-8}
\vsp 15& & & & &1&22& 441  \\
\cline{1-8}
\vsp 16& & & & & & 7& 345  \\
\cline{1-8}
\vsp 17& & & & & & 1& 192 \\
\cline{1-8}
\vsp 18& & & & & &  & \it 75 \\
\cline{1-8}
\vsp 19& & & & & &  &  \it 12  \\
\cline{1-8}
\vsp 20& & & & & &  &  \it 3  \\
\cline{1-8}
\vsp 21& & & & & &  &   \it 3  \\
\cline{1-8}
\vsp 22& & & & & &  &    \\
\cline{1-8}
\vsp 23& & & & & &  &    \\
\cline{1-8}
\vsp 24& & & & & &  &   1 \\
\cline{1-8}
\cline{1-8}
\vsp Subtotal& 1& 1& 3& 8& {29}& {240}& {2335}\\
\cline{1-8}
\end{tabular}
\\ \medskip
\caption{Number of prime knots with given arc index and minimal crossing number}
\end{table}
}

The authors acknowledge the effort of Alexander Stoimenow in eliminating duplicates from the table among the knots with more than 16 crossings.
In~\cite{N1999}, Nutt identified all knots up to arc index 9. In~\cite{B2002}, Beltrami determined arc index for prime knots up to ten crossings. In~\cite{Jin2006}, Jin et al.\ identified all prime knots up to arc index 10. Ng determined arc index for prime knots up to eleven crossings~\cite{Ng2006}. Matsuda determined the arc index for torus knots~\cite{H2006}. The \emph{Table of Knot Invariants\/}~\cite{knotinfo} gives the arc index for prime knots up to 12 crossings.

In the next section, we list the prime knots whose arc index is less than or equal to 11 together with corresponding minimal arc presentations in the form of projections of the grid diagrams. The names of the knots used are the Dowker-Thistlethwaite names\footnote{The Knotscape names~\cite{knotscape,knotinfo}} (or DT names) up to 16 crossings. Up to 10 crossings, we also put the classical names\footnote{Adjusted from names in the book of Rolfsen after the discovery of the Perko pair~\cite{rolfsen,knotinfo}} in brackets.
For knots with more than 16 crossings, we named them $17n_1, 17n_2,\ldots,18n_1,\ldots$, as they appear in the list. We also put the Dokwer-Thistlethwaite notation (or DT notation) for the knots with more than 16 crossings.

\begin{figure}[h]
\centering
\setlength{\unitlength}{0.11mm}
\myfbox{
\begin{picture}(240,240)(-20,-20)
\thicklines
\put(  0, 80){\line( 0,1){120}}
\put(  0,200){\line( 1,0){160}}
\put(  0, 80){\line( 1,0){ 30}}
\put( 40, 40){\line( 0,1){ 80}}
\put( 40,120){\line( 1,0){ 30}}
\put( 40, 40){\line( 1,0){ 70}}
\put( 80, 80){\line( 0,1){ 80}}
\put( 80,160){\line( 1,0){ 70}}
\put( 80, 80){\line(-1,0){ 30}}
\put(120,  0){\line( 0,1){120}}
\put(120,120){\line(-1,0){ 30}}
\put(120,  0){\line( 1,0){ 80}}
\put(160, 40){\line( 0,1){160}}
\put(160, 40){\line(-1,0){ 30}}
\put(200,  0){\line( 0,1){160}}
\put(200,160){\line(-1,0){ 30}}
\end{picture}
}
\qquad
\myfbox{
\begin{picture}(240,240)(-20,-20)
\qbezier[40](-20,  0)(100,  0)(220,  0)
\qbezier[40](-20, 40)(100, 40)(220, 40)
\qbezier[40](-20, 80)(100, 80)(220, 80)
\qbezier[40](-20,120)(100,120)(220,120)
\qbezier[40](-20,160)(100,160)(220,160)
\qbezier[40](-20,200)(100,200)(220,200)
\qbezier[40](  0,-20)(  0,100)(  0,220)
\qbezier[40]( 40,-20)( 40,100)( 40,220)
\qbezier[40]( 80,-20)( 80,100)( 80,220)
\qbezier[40](120,-20)(120,100)(120,220)
\qbezier[40](160,-20)(160,100)(160,220)
\qbezier[40](200,-20)(200,100)(200,220)
\scriptsize
{
\put( -15,204){$1$}\put( 25,204){$0$}\put( 65,204){$0$}\put( 105,204){$0$}\put( 145,204){$1$}\put( 185,204){$0$}
\put( -15,164){$0$}\put( 25,164){$0$}\put( 65,164){$1$}\put( 105,164){$0$}\put( 145,164){$0$}\put( 185,164){$1$}
\put( -15,124){$0$}\put( 25,124){$1$}\put( 65,124){$0$}\put( 105,124){$1$}\put( 145,124){$0$}\put( 185,124){$0$}
\put( -15, 84){$1$}\put( 25, 84){$0$}\put( 65, 84){$1$}\put( 105, 84){$0$}\put( 145, 84){$0$}\put( 185, 84){$0$}
\put( -15, 44){$0$}\put( 25, 44){$1$}\put( 65, 44){$0$}\put( 105, 44){$0$}\put( 145, 44){$1$}\put( 185, 44){$0$}
\put( -15,  4){$0$}\put( 25,  4){$0$}\put( 65,  4){$0$}\put( 105,  4){$1$}\put( 145,  4){$0$}\put( 185,  4){$1$}
}%
{
\thicklines
\put(  0, 80){\line( 0,1){120}}
\put(  0,200){\line( 1,0){160}}
\put(  0, 80){\line( 1,0){ 30}}
\put( 40, 40){\line( 0,1){ 80}}
\put( 40,120){\line( 1,0){ 30}}
\put( 40, 40){\line( 1,0){ 70}}
\put( 80, 80){\line( 0,1){ 80}}
\put( 80,160){\line( 1,0){ 70}}
\put( 80, 80){\line(-1,0){ 30}}
\put(120,  0){\line( 0,1){120}}
\put(120,120){\line(-1,0){ 30}}
\put(120,  0){\line( 1,0){ 80}}
\put(160, 40){\line( 0,1){160}}
\put(160, 40){\line(-1,0){ 30}}
\put(200,  0){\line( 0,1){160}}
\put(200,160){\line(-1,0){ 30}}
}
\end{picture}
}
\qquad
\myfbox{
\begin{picture}(240,240)(-20,-20)
\scriptsize
{
\put( -15,204){$1$}\put( 25,204){$0$}\put( 65,204){$0$}\put( 105,204){$0$}\put( 145,204){$1$}\put( 185,204){$0$}
\put( -15,164){$0$}\put( 25,164){$0$}\put( 65,164){$1$}\put( 105,164){$0$}\put( 145,164){$0$}\put( 185,164){$1$}
\put( -15,124){$0$}\put( 25,124){$1$}\put( 65,124){$0$}\put( 105,124){$1$}\put( 145,124){$0$}\put( 185,124){$0$}
\put( -15, 84){$1$}\put( 25, 84){$0$}\put( 65, 84){$1$}\put( 105, 84){$0$}\put( 145, 84){$0$}\put( 185, 84){$0$}
\put( -15, 44){$0$}\put( 25, 44){$1$}\put( 65, 44){$0$}\put( 105, 44){$0$}\put( 145, 44){$1$}\put( 185, 44){$0$}
\put( -15,  4){$0$}\put( 25,  4){$0$}\put( 65,  4){$0$}\put( 105,  4){$1$}\put( 145,  4){$0$}\put( 185,  4){$1$}
}%
\end{picture}
}
\caption{Correspondence of grid diagram and Cromwell matrix}\label{fig:grid2matrix}
\end{figure}
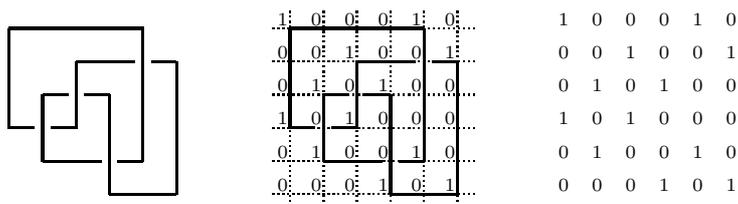

The \emph{norm\/} of an arc presentation with $n$ arcs is the $n^2$ digit binary number which can be read from the corresponding Cromwell matrix.
The norm of the arc presentation of Figure~\ref{fig:grid2matrix} equals $100010\ 001001\ 010100\ 101000\ 010010\ 000101_2$. The diagrams in the next section are the ones with the maximal norm for each knot type. The knots with arc index 11 and more than 16 crossings are listed in the norm-decreasing order of the given arc presentations for each crossing number.
\clearpage
\renewcommand{\thefootnote}{\fnsymbol{footnote}}
\section{A tabulation of primes knots up to arc index 11 and their minimal arc presentations$^\dagger$}
\footnotetext[2]{Minimal arc presentations are given by the projections of minimal grid diagrams in which the vertical strands cross over the horizontal strands. Knot types up to 16 crossings were determined by using Knotscape. Knot types of 17 crossings were estimated by Knotscape and then confirmed by distinguishing them from knots of smaller crossings.}

\setlength{\unitlength}{6.2pt}
\newcommand{\myspace}{\ }
\newcommand{\mythicklines}{}
\newcommand{\mythinlines}{}
\newcommand{\rolfsen}[2]{{[$#1_{#2}$]}}

\noindent
{\bf Arc Index 5}\nopagebreak

\noindent
 \
\dt{4 & 14 & -20 & 22 \cr -30 & 24 & 2 & -26 \cr -32 & -34 & -28 & 8 \cr -16 & 12 & -6 & -10 \cr -18 &  &  & } 

\medskip\setcounter{cr}{18}\setcounter{knumber}{0}

\noindent{\bf Arc Index 11} (18 crossings$^\dagger$)\nopagebreak
\footnotetext[2]{The crossing number for these knots were estimated by using Knotscape. It is certain that the minimal crossing number for these knots is bigger than 16.}

\nopagebreak\noindent
%
 \ \dt{} 

\medskip\setcounter{cr}{19}\setcounter{knumber}{0}

\noindent{\bf Arc Index 11} (19 crossings$^\dagger$)\nopagebreak
\footnotetext[2]{The crossing number for these knots were estimated by using Knotscape. It is certain that the minimal crossing number for these knots is bigger than 16.}

\noindent
%
 \
\dt{4 & 14 & -22 & 20 \cr -28 & 18 & 2 & -24 \cr 36 & 10 & -30 & -38 \cr 34 & -12 & -8 & 6 \cr 16 & 26 & 32 & } 

\medskip\setcounter{cr}{20}\setcounter{knumber}{0}

\noindent{\bf Arc Index 11} (20 crossings$^\dagger$)\nopagebreak
\footnotetext[2]{The crossing number for these knots were estimated by using Knotscape. It is certain that the minimal crossing number for these knots is bigger than 16.}

\noindent
%
 \ \dt{} 

\medskip\setcounter{cr}{21}\setcounter{knumber}{0}

\noindent{\bf Arc Index 11} (21 crossings$^\dagger$)\nopagebreak
\footnotetext[2]{As $21n_1$ is the torus knot of type (4,7), its minimal crossing number is 21. The crossing number for the other two knots were estimated by using Knotscape. It is certain that the minimal crossing number for them is bigger than 16.}

\nopagebreak\noindent
%
 \ \dt{} 

\medskip\setcounter{cr}{24}\setcounter{knumber}{0}

\noindent{\bf Arc Index 11} (24 crossings$^\ddagger$)\nopagebreak
\footnotetext[3]{As $24n_1$ is the torus knot of type (5,6), its minimal crossing number is 24.}

\noindent
%
 \
\phantom{\dt{}}


\clearpage
\begin{thebibliography}{MOT}
\bibitem{BP2000}
            Yongju Bae and Chan-Young Park,
            \emph{An upper bound of arc index of links},
            Math. Proc. Camb. Phil. Soc. \textbf{129} (2000) 491--500.
\bibitem{BG2006}
            J. A. Baldwin and W. D. Gillam,
            \emph{Computations of Heegaard.Floer knot homology},
            arXiv: math/0610167.
\bibitem{B2002}
            Elisabeta Beltrami, \emph{Arc index of non-alternating links},
            Knots 2000 Korea vol.1, J. Knot Theory Ramifications. \textbf{11}(3) (2002) 431--444.
\bibitem{C1995}
            Peter R. Cromwell, \emph{Embedding knots and links in an open book I: Basic properties},
            Topology Appl. \textbf{64} (1995) 37--58.
\bibitem{CN1996}
            Peter R. Cromwell and Ian J. Nutt,
            \emph{Embedding knots and links in an open book II: Bounds on arc index},
            Math. Proc. Camb. Phil. Soc. \textbf{119} (1996), 309--319.
\bibitem{HTW1998}
            Jim Hoste, Morwen Thistlethwaite and Jeff Weeks,
            \emph{The first 1,701,936 knots},
            Math. Intelligencer \textbf{20}(4) (1998) 33--48.
\bibitem{Jin2006}
            Gyo Taek Jin, Hun Kim, Gye-Seon Lee, Jae Ho Gong, Hyuntae Kim, Hyunwoo Kim and Seul Ah Oh,
            \emph{Prime knots with arc index up to 10},
            Intelligence of Low Dimensional Topology 2006,
            Series on Knots and Everything Book vol. 40, World Scientific Publishing Co., 65--74, 2006.
\bibitem{Jin2007}
            Gyo Taek Jin and Wang Keun Park,
            \emph{Prime knots with arc index up to 11 and an upper bound of arc index for non-alternating knots},
            J. Knot Theory Ramifications, to appear.
\bibitem{H2006}
            Hiroshi Matsuda,
            \emph{Links in an open book decomposition and in the standard contact structure},
            Proc. Amer. Math. Soc. \textbf{134}(12) (2006) 3697--3702 (electronic).
\bibitem{Ng2006}
            Lenhard Ng, \emph{On arc index and maximal Thurston-Bennequin number},
            arXiv: math/0612356
\bibitem{N1997}
            Ian J. Nutt, \emph{Arc index and Kauffman polynomial},
            J. Knot Theory Ramifications. \textbf{6}(1) (1997) 61--77.
\bibitem{N1999}
            Ian J. Nutt,
            \emph{Embedding knots and links in an open book III: On the braid index of satellite links},
            Math. Proc. Camb. Phil. Soc. \textbf{126} (1999) 77--98.
\bibitem{rolfsen}
            Dale Rolfsen, \emph{Knots and Links}, AMS Chelsea Publishing, 2003
\bibitem{knotscape}
            Knotscape, \texttt{http://www.math.utk.edu/$\sim$morwen/knotscape.html}
\bibitem{knotinfo}
            Table of Knot Invariants, \texttt{http://www.indiana.edu/$\sim$knotinfo/}

\end{thebibliography}
\end{document}